\documentclass[12pt]{article} 

\usepackage{amsmath}  
\usepackage{amsfonts} 
\usepackage{graphicx} 
\usepackage{url}
\usepackage[bottom]{footmisc}
\usepackage{amssymb}
\usepackage{booktabs}
\usepackage[labelformat=simple]{subcaption}
\usepackage{textcomp}
\usepackage{mathrsfs}
\usepackage[shortlabels]{enumitem}
\usepackage{upgreek}
\usepackage{fullpage}
\usepackage[title]{appendix}
\usepackage[margin=1in]{geometry}
\usepackage{color}
\usepackage{times}
\usepackage{pgf,tikz}
\usepackage[raggedright]{titlesec}
\usepackage{multirow}
\usepackage{multicol}
\usepackage{array}
\usepackage{setspace}
\usepackage[format=plain,labelsep=period]{caption}
\usepackage{cite}
\usepackage{esint}
\usepackage{listings}

\titleformat{\section}{\bfseries}{\thesection}{1em}{}
\titleformat{\subsection}{\itshape}{\thesubsection}{1em}{}
\titleformat{\subsubsection}{}{\thesubsubsection}{1em}{}

\newcolumntype{L}[1]{>{\raggedright\let\newline\\\arraybackslash\hspace{0pt}}m{#1}}
\newcolumntype{C}[1]{>{\centering\let\newline\\\arraybackslash\hspace{0pt}}m{#1}}
\newcolumntype{R}[1]{>{\raggedleft\let\newline\\\arraybackslash\hspace{0pt}}m{#1}}

\begin{document}

\begin{table}[!t]
\centering
\begin{tabular}{c}
{\large \textbf{Unconditionally stable symplectic integrators for the Navier-Stokes}} \\
{\large \textbf{equations and other dissipative systems}} \\ \\
\end{tabular}
\end{table}

\begin{table}[!h]
\centering
\begin{tabular}{cccc}
Sutthikiat Sungkeetanon & Joseph S. Gaglione & Robert L. Chapman & Tyler M. Kelly  \\
Howard A. Cushman & Blakeley H. Odom & Bryan MacGavin & Gafar A. Elamin \\
Nathan J. Washuta & Jonathan E. Crosmer & Adam C. DeVoria & John W. Sanders$^{*}$
\end{tabular}
\end{table}

\begin{table}[!h]
\centering
\begin{tabular}{c}
{Department of Mechanical Engineering} \\
{The Citadel, The Military College of South Carolina} \\
{171 Moultrie St, Charleston, SC 29409} \\ \\
{$^{*}$Corresponding author: jsande12@citadel.edu} \\
\end{tabular}
\end{table}

\section*{Abstract}

Symplectic integrators offer vastly superior performance over traditional numerical techniques for conservative dynamical systems, but their application to \emph{dissipative} systems is inherently difficult due to dissipative systems' lack of symplectic structure. Leveraging the intrinsic variational structure of higher-order dynamics, this paper presents a general technique for applying existing symplectic integration schemes to dissipative systems, with particular emphasis on viscous fluids modeled by the Navier-Stokes equations. Two very simple such schemes are developed here. Not only are these schemes unconditionally stable for dissipative systems, they also outperform traditional methods with a similar degree of complexity in terms of accuracy for a given time step. For example, in the case of viscous flow between two infinite, flat plates, one of the schemes developed here is found to outperform both the implicit Euler method and the explicit fourth-order Runge-Kutta method in predicting the velocity profile. To the authors' knowledge, this is the very first time that a symplectic integration scheme has been applied successfully to the Navier-Stokes equations. We interpret the present success as direct empirical validation of the canonical Hamiltonian formulation of the Navier-Stokes problem recently published by Sanders~\emph{et al.} More sophisticated symplectic integration schemes are expected to exhibit even greater performance. It is hoped that these results will lead to improved numerical methods in computational fluid dynamics. 



\section{Introduction and literature review}

Symplectic integrators are a class of incredibly powerful numerical integration schemes that offer vastly superior performance over traditional numerical methods for conservative, Hamiltonian systems~\cite{Hairer2006}. Built on the seminal work of several pioneering researchers including Maeda~\cite{Maeda1980,Maeda1982}, Veselov~\cite{Veselov1988,Veselov1991}, Suris~\cite{Suris1988,Suris1990,Suris1996,Suris2003}, and MacKay~\cite{MacKay1992}, these new and improved algorithms experienced a remarkable surge of interest in recent decades with the prolific work of Marsden and coworkers~\cite{Wendlandt1997,Kane2000,Marsden2001,Lew2003,Lew2004,Lew2004a,Kale2007,OberBloebaum2011,Lew2016}. As the name suggests, symplectic integrators owe their robust performance to the preservation of the symplectic structure of Hamiltonian dynamical systems~\cite{Suris1988}. Consequently, relatively little work has been done by way of applying symplectic integrators to non-Hamiltonian systems, which by their very nature \emph{lack} symplectic structure. Indeed, such endeavors must usually rely on \emph{ad hoc} pseudo-variational principles such as the d'Alembert-Lagrange principle~\cite{Kraus2015,Kraus2021}. One notable exception is recent work by Limebeer and coworkers~\cite{Limebeer2020,Limebeer2023}, which utilizes a ``transmission line'' composed of several springs and inerters in series to model a heat bath for both linear~\cite{Limebeer2020} and nonlinear~\cite{Limebeer2023} dissipation. 

The present work takes an altogether different approach to applying symplectic integrators to dissipative systems which leverages the intrinsic variational structure of higher-order dynamics~\cite{Sanders2023a,Sanders2023b}. It was recently discovered that the principle of least squares---when the sum of squares is averaged over an arbitrary time period---can transform any non-Hamiltonian system into a Hamiltonian one by doubling the order of the equation~\cite{Sanders2021,Sanders2022,Sanders2023a,Sanders2023,Sanders2023b,Sandersetal2024,SandersBecker2024,SandersBecker2024a,Becker2024}. Remarkably, even if the dynamical system lacks symplectic structure in its original $(2\mathcal{N}+1)$-dimensional phase space, it nevertheless {does} possess symplectic structure in the $(4\mathcal{N}+1)$-dimensional phase space of the higher-order formulation. As we will see presently, this approach leads to a completely general framework whereby existing symplectic integration schemes (originally intended for conservative systems) may be applied to any arbitrary dissipative system---including, most notably, viscous fluids modeled by the Navier-Stokes equations~\cite{Stokes1845,LemarieRieusset2018,Sandersetal2024}.

One popular method for constructing symplectic integrators is based on a discretized version of Hamilton's principle~\cite{Hairer2006}. Following Chapter VI.6 of the text by Hairer~\emph{et~al.}~\cite{Hairer2006}, suppose we have a dynamical system with generalized coordinates $q(t)$ (and conjugate momenta $p(t)$) whose dynamics are described by a Lagrangian $L(t,q,\dot{q})$. Upon discretizing time $t\mapsto\{t_{n}; n=0,1,2,...,N\}$, the action integral over any time interval $t\in(0,\tau)$ may be approximated by a discrete sum
\begin{equation}
\mathcal{S}\equiv\int_{0}^{\tau}L(t,q,\dot{q})\text{d}t\approx\sum_{n=0}^{N-1}L_{h}(t_{n},t_{n+1},q_{n},q_{n+1})\equiv\mathcal{S}_{h},
\end{equation}
where $q_{n}\equiv q(t_{n})$ and the \emph{discrete Lagrangian} $L_{h}$ is given by an approximation to the integral
\begin{equation}\label{eq:}
L_{h}(t_{n},t_{n+1},q_{n},q_{n+1})\approx\int_{t_{n}}^{t_{n+1}}L(t,q,\dot{q})\text{d}t.
\end{equation}
The condition that $\mathcal{S}_{h}$ be stationary with respect to each of the discrete coordinates $q_{n}$ leads to the \emph{discrete Euler-Lagrange equations}
\begin{equation}\label{eq:discreteEL}
\frac{\partial L_{h}}{\partial y}(t_{n-1},t_{n},q_{n-1},q_{n})+\frac{\partial L_{h}}{\partial x}(t_{n},t_{n+1},q_{n},q_{n+1})=0, \quad n=1,...,N-1,
\end{equation}
where $x$ and $y$ refer to the third and fourth arguments, respectively, of $L_{h}(t_{n},t_{n+1},x,y)$.

In analogy to Hamilton's principal function, denote by $S_{h}$ the result upon evaluating the discrete action $\mathcal{S}_{h}$ for discrete coordinates $q_{n}$ satisfying the discrete Euler-Lagrange equations:
\begin{equation}
S_{h}(q_{0},q_{N})\equiv\sum_{n=0}^{N-1}L_{h}(t_{n},t_{n+1},q_{n},q_{n+1}).
\end{equation}
The \emph{discrete momenta} $p_{n}$ are defined as
\begin{equation}\label{eq:pold}
p_{n}\equiv-\frac{\partial L_{h}}{\partial q_{n}}(t_{n},t_{n+1},q_{n},q_{n+1}),
\end{equation}
so that, using the discrete Euler-Lagrange equations~\eqref{eq:discreteEL},
\begin{equation}\label{eq:pnew}
p_{n+1}=+\frac{\partial L_{h}}{\partial q_{n+1}}(t_{n},t_{n+1},q_{n},q_{n+1}),
\end{equation}
and
\begin{equation}
\text{d}S_{h}=p_{N}\text{d}q_{N}-p_{0}\text{d}q_{0},
\end{equation}
once again in direct analogy to the classical Hamilton-Jacobi theory. It follows that the $p_{n}$ are indeed discrete approximations to the conjugate momenta $p(t)$.

At each time step $t_{n}$, one knows $q_{n}$ and $p_{n}$. In general, one solves~\eqref{eq:pold} for $q_{n+1}$ and then evaluates $p_{n+1}$ directly using~\eqref{eq:pnew}, resulting in the symplectic map~\cite{Maeda1980,Suris1990,Veselov1991,MacKay1992}
\begin{equation}\label{eq:}
(q_{n},p_{n})\xrightarrow{\text{\eqref{eq:pold}}}q_{n+1}\xrightarrow{\text{\eqref{eq:pnew}}}p_{n+1}.
\end{equation}
Distinct symplectic integrators are obtained by selecting particular approximations for the discrete Lagrangian $L_{h}$ and the velocities $\dot{q}$. For instance, one of the simplest schemes, due to MacKay~\cite{MacKay1992}, employs the trapezoidal rule for $L_{h}$ with a constant time increment $t_{n+1}-t_{n}=h$, along with a forward finite difference approximation for $\dot{q}$. Thus, with MacKay's method~\cite{MacKay1992}, the discrete Lagrangian takes the form
\begin{equation}\label{eq:}
L_{h}(t_{n},t_{n+1},q_{n},q_{n+1})=\frac{h}{2}L\left(t_{n},q_{n},\frac{q_{n+1}-q_{n}}{h}\right)+\frac{h}{2}L\left(t_{n+1},q_{n+1},\frac{q_{n+1}-q_{n}}{h}\right).
\end{equation}
As noted by Hairer~\emph{et~al.}~\cite{Hairer2003}, when applied to the classical Lagrangian $L=T-V$ (where $T$ and $V$ are the kinetic and potential energies, respectively), MacKay's method~\cite{MacKay1992} reproduces the well known St\"{o}rmer-Verlet method~\cite{Hairer2003,Hairer2006}. Other symplectic integrators of note include the method of Wendlandt~\&~Marsden~\cite{Wendlandt1997}, which employs the midpoint rule rather than the trapezoidal rule, and the method of Marsden~\&~West~\cite{Marsden2001}, which leads to a symplectic, partitioned Runge-Kutta method.

Symplectic integrators as such can only be applied directly to dynamical systems with symplectic structure. In the context of mechanics, ``symplectic'' is synonymous with ``conservative'' and ``Hamiltonian.'' Because non-Hamiltonian systems lack symplectic structure, it is not possible to apply symplectic integrators directly to non-Hamiltonian systems~\cite{Kraus2015,Kraus2021,Limebeer2020,Limebeer2023}.

The remainder of this paper is organized as follows. Section~\ref{sec:general} describes a general technique for applying any existing symplectic integrator in the literature to generic dynamical systems---whether conservative or dissipative---via a higher-order formulation. Section~\ref{sec:stability} develops two specific such schemes based on the method of MacKay~\cite{MacKay1992}. It is shown that both schemes are unconditionally stable and generally outperform the implicit Euler method. Section~\ref{sec:quadraticdrag} applies the two aforementioned schemes to the nonlinear problem of quadratic drag in one dimension, with similar results. Section~\ref{sec:poiseuille} applies these schemes to the Navier-Stokes equations for the case of developing viscous flow between two infinite, fixed plates. Notably, one of the two schemes developed here is found to outperform \emph{both} the implicit Euler method and the explicit fourth-order Runge-Kutta method. Finally, Section~\ref{sec:conclusion} concludes the paper with a brief summary and some closing remarks. 

\section{Symplectic integrators for general dissipative systems}\label{sec:general}

Consider the generic first-order problem
\begin{equation}\label{eq:genericfirstorder}
\dot{v}=f(t,v), \quad v(0)=v_{0},
\end{equation}
where $f(t,v)$ is an arbitrary function of $t$ and $v$. In the specific context of mechanics, if the forces on a dissipative system depend only on the system's velocities and possibly time, but not on the system's configuration, then we identify $v(t)$ with the {velocities} of said system. Notably, such is the case for the Navier-Stokes equations~\cite{Stokes1845,LemarieRieusset2018}. More generally, because any system of differential equations can always be written as a first-order system via the introduction of new variables,~\eqref{eq:genericfirstorder} can be considered representative of arbitrary dynamical systems (including dissipative ones). 

Note that~\eqref{eq:genericfirstorder} does \emph{not} possess symplectic structure, since it is a known fact that first-order time derivatives cannot be obtained by the variation of a first-order Lagrangian without introducing non-physical terms into the governing equation~\cite{SandersBecker2024a}. As such, we cannot apply existing symplectic integrators \emph{directly} to~\eqref{eq:genericfirstorder}. Nevertheless, we {can} apply symplectic integrators to a mathematically equivalent higher-order formulation by leveraging the intrinsic variational structure of higher-order dynamics~\cite{Sanders2023a,Sanders2023b}.

Following Sanders~\cite{Sanders2023a,Sanders2023b}, we observe that the solution $v(t)$ to~\eqref{eq:genericfirstorder} corresponds to a local minimum of the following action:
\begin{equation}\label{eq:}
\mathcal{S}^{*}[v]=\int_{0}^{\tau}\tfrac{1}{2}\mathcal{R}^{2}\text{d}t=\int_{0}^{\tau}\tfrac{1}{2}\left(\dot{v}-f(t,v)\right)^{2}\text{d}t,
\end{equation}
where $\mathcal{R}\equiv\dot{v}-f(t,v)$ is the residual. Indeed, the reader can easily verify that $\delta \mathcal{S}^{*}=\int\mathcal{R}\delta\mathcal{R}\text{d}t=0$ for the actual motion (for which the residual vanishes), and $\delta^{2}\mathcal{S}^{*}=\int(\delta\mathcal{R})^{2}\text{d}t>0$. This is the time-averaged principle of least squares~\cite{Sanders2021,Sanders2022,Sanders2023a,Sanders2023,Sanders2023b,Sandersetal2024,SandersBecker2024,SandersBecker2024a,Becker2024}.  Varying $v(t)$, we have that
\begin{equation}\label{eq:}
\delta\mathcal{S}^{*}=\int_{0}^{\tau}\left(\dot{v}-f(t,v)\right)\left(\delta\dot{v}-f_{v}(t,v)\delta v\right)\text{d}t,
\end{equation}
where $f_{v}\equiv\partial f/\partial v$. Upon integration by parts, we find that
\begin{equation}\label{eq:EL}
\delta\mathcal{S}^{*}=\int_{0}^{\tau}\left(-\ddot{v}+{f}_{t}(t,v)+f_{v}(t,v)f(t,v)\right)\delta v\text{d}t+\left[\left(\dot{v}-f(t,v)\right)\delta v\right]_{0}^{\tau},
\end{equation}
where $f_{t}\equiv\partial f/\partial t$. The Euler-Lagrange equation is obtained by setting $\delta\mathcal{S}^{*}=0$ for arbitrary $\delta v(t)$, yielding the \emph{second-order} equation
\begin{equation}\label{eq:secondorderEOM}
\ddot{v}={f}_{t}(t,v)+f_{v}(t,v)f(t,v).
\end{equation}
With initial conditions $v(0)=v_{0}$ and $\dot{v}(0)=f(0,v_{0})$, the solution of~\eqref{eq:secondorderEOM} will coincide with the solution of~\eqref{eq:genericfirstorder}. Furthermore, note that the ``momentum'' conjugate to $v(t)$---which is read from the boundary term in~\eqref{eq:EL}---coincides precisely with the residual: $p(t)\equiv\dot{v}-f(t,v)\equiv\mathcal{R}$. As we will see presently, this observation is crucial and holds the key to ensuring unconditional stability.

Even though the first-order problem~\eqref{eq:genericfirstorder} lacks symplectic structure, nevertheless the mathematically equivalent second order problem~\eqref{eq:secondorderEOM} {is} symplectic within its higher-dimensional phase space $(v,p,t)$. Indeed, Liouville's theorem is automatically satisfied, since the motion is confined to the plane defined by $p=0$, so that the phase-space volume, being always zero, is conserved~\cite{Sandersetal2024}. In this way, it is possible to apply existing symplectic integrators to the dissipative system~\eqref{eq:genericfirstorder} via the equivalent second-order formulation~\eqref{eq:secondorderEOM}. The associated Lagrangian is
\begin{equation}\label{eq:L}
L(t,v,\dot{v})\equiv\tfrac{1}{2}\left(\dot{v}-f(t,v)\right)^{2}.
\end{equation}
To this Lagrangian may be applied any existing symplectic integrator in the literature. In what follows, we will focus on the particular method due to MacKay~\cite{MacKay1992}, in order to illustrate just how powerful even the simplest of symplectic integrators can be when applied to dissipative systems. We expect more sophisticated schemes, such as that of Marsden~\&~West~\cite{Marsden2001}, to perform even better.

\section{Test problem and stability analysis}\label{sec:stability}

Let us consider the test problem for stability analysis, which, after non-dimensionalization, reads
\begin{equation}\label{eq:test}
\dot{v}=-v, \quad v(0)=1.
\end{equation}
Note that this can be interpreted as the problem of a lumped mass moving in a viscous medium with linear drag. This is an inherently non-Hamiltonian problem that fundamentally lacks symplectic structure, which is why it is not possible to apply symplectic integrators directly to this problem.

The residual is $\mathcal{R}\equiv\dot{v}+v$, and so the Lagrangian of the mathematically equivalent second-order problem is
\begin{equation}
L(v,\dot{v})\equiv\frac{1}{2}\dot{v}^{2}+v\dot{v}+\frac{1}{2}v^{2}.
\end{equation}
Again, at this point we may select any existing symplectic integration scheme. For the sake of simplicity, here we will employ the method of MacKay~\cite{MacKay1992}. Recall that MacKay's method~\cite{MacKay1992} employs the trapezoidal rule for the time integral with constant time increment $h$, along with a forward finite difference approximation for $\dot{v}$. Thus, the discrete Lagrangian is given by
\begin{equation}\label{eq:discreteLMacKay}
L_{h}(v_{n},v_{n+1})=\frac{h}{2}L\left(v_{n},\frac{v_{n+1}-v_{n}}{h}\right)+\frac{h}{2}L\left(v_{n+1},\frac{v_{n+1}-v_{n}}{h}\right).
\end{equation}
The discrete momenta are found to be
\begin{equation}\label{eq:poldMacKaytest}
p_{n}=-\frac{\partial L_{h}}{\partial v_{n}}(v_{n},v_{n+1})=\frac{1}{h}v_{n+1}-\left(\frac{1}{h}-1+\frac{h}{2}\right)v_{n},
\end{equation}
\begin{equation}\label{eq:pnewMacKaytest}
p_{n+1}=+\frac{\partial L_{h}}{\partial v_{n+1}}(v_{n},v_{n+1})=\left(\frac{1}{h}+1+\frac{h}{2}\right)v_{n+1}-\frac{1}{h}v_{n}.
\end{equation}
A ``conventional'' application of MacKay's method~\cite{MacKay1992}---that is, at each time step solving~\eqref{eq:poldMacKaytest} for $v_{n+1}$ in terms of the known $v_{n}$ and $p_{n}$, and then using~\eqref{eq:pnewMacKaytest} to evaluate $p_{n+1}$ for the next time step---yields an unstable numerical solution. In order to ensure stability, we will take advantage of the fact that, within the context of the second-order formulation, we know \emph{a priori} the exact value of the conjugate momentum. 

Recall that, in the second-order formulation, the momentum conjugate to the velocity coincides with the residual, which vanishes for the actual motion. Thus, $p(t)=0$ for the actual motion. Now, for a fixed $v_{n}$, the only way to have \emph{both} $p_{n}=0$ \emph{and} $p_{n+1}=0$ is to make $h=0$. It would seem, then, that we must choose \emph{either} $p_{n}=0$ \emph{or} $p_{n+1}=0$.

When we set $p_{n}=0$, we obtain
\begin{equation}\label{eq:}
v_{n+1}=\left({1-h+\tfrac{1}{2}h^{2}}\right)v_{n},
\end{equation}
which we recognize as the explicit Euler method with a quadratic correction term $\tfrac{1}{2}h^{2}$. Likewise, when we set $p_{n+1}=0$, we obtain
\begin{equation}\label{eq:}
v_{n+1}=\left(\frac{1}{1+h+\tfrac{1}{2}h^{2}}\right)v_{n},
\end{equation}
which we recognize as the \emph{implicit} Euler method with quadratic correction. Notably, the implicit Euler method is unconditionally stable. Indeed, since
\begin{equation}\label{eq:}
\left|\frac{1}{1+h+\tfrac{1}{2}h^{2}}\right|<1
\end{equation}
for all $h>0$, we have established that a symplectic integrator obtained by applying MacKay's method~\cite{MacKay1992} to the second-order Lagrangian~\eqref{eq:L} and setting $p_{n+1}=0$ will be unconditionally stable, even for dissipative systems. We note that, as $h\rightarrow\infty$, this method gives the correct steady-state value $v_{n+1}=0$ for the test problem.

Of course, there are alternatives to setting either $p_{n}=0$ or $p_{n+1}=0$. For instance, motivated by a desire to ensure that the conjugate momentum is conserved, we could insist that the discrete momentum is {constant}, so that the \emph{difference} between $p_{n}$ and $p_{n+1}$ vanishes:
\begin{equation}
p_{n+1}-p_{n}=\left(1+\frac{h}{2}\right)v_{n+1}-\left(1-\frac{h}{2}\right)v_{n}=0.
\end{equation} 
This yields
\begin{equation}
v_{n+1}=\left(\frac{1-\tfrac{h}{2}}{1+\tfrac{h}{2}}\right)v_{n},
\end{equation}
which is also unconditionally stable for $h>0$. In this case, as $h\rightarrow\infty$, we have that $v_{n+1}=-v_{n}$, which for finite $h$ would exhibit oscillation about the steady-state value $v=0$.

Following this line of reasoning, we could insist that any linear combination of $p_{n}$ and $p_{n+1}$ vanishes:
\begin{equation}
\alpha p_{n}+\beta p_{n+1}=-\alpha\left(\frac{1+\beta/\alpha}{h}-1+\frac{h}{2}\right)v_{n}+\beta\left(\frac{1+\alpha/\beta}{h}+1+\frac{h}{2}\right)v_{n+1}=0,
\end{equation}
which yields
\begin{equation}
v_{n+1}=\left(\frac{\alpha+\beta-\alpha h+\tfrac{1}{2}\alpha h^{2}}{\alpha+\beta+\beta h+\tfrac{1}{2}\beta h^{2}}\right)v_{n}.
\end{equation}
Different choices of $(\alpha,\beta)$ will in general yield different schemes with different levels of performance. We leave a detailed investigation into the optimal choice of $(\alpha,\beta)$ for future work. Henceforth, we will refer to the method characterized by $(\alpha,\beta)=(0,1)$ as ``Method~I'' (this corresponds to setting $p_{n+1}=0$) and the method characterized by $(\alpha,\beta)=(-1,1)$ as ``Method~II'' (this corresponds to setting $p_{n+1}-p_{n}=0$). Not only are these two methods unconditionally stable, we also expect them to outperform the implicit Euler method in terms of both accuracy and rate of convergence, since the present methods will in general possess higher-order correction terms that the standard implicit Euler method lacks. \emph{Indeed, the time-averaged principle of least squares completely automates the inclusion of these higher-order correction terms}.

\begin{figure}[!p]
\centering
\begin{subfigure}{0.45\textwidth}
\centering
\includegraphics[width=\textwidth]{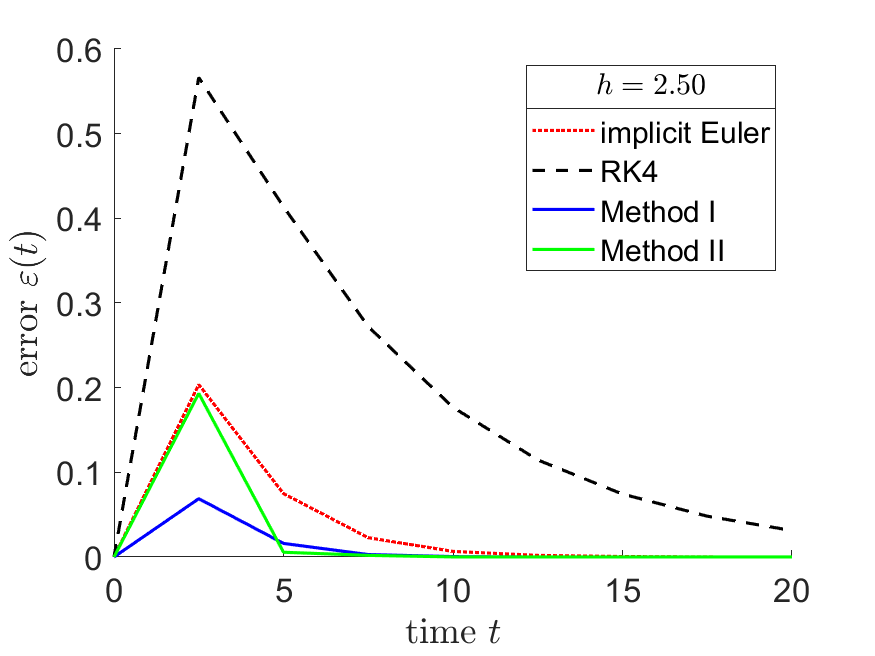}
\caption{}
\label{fig:fig1a}
\end{subfigure}
\begin{subfigure}{0.45\textwidth}
\centering
\includegraphics[width=\textwidth]{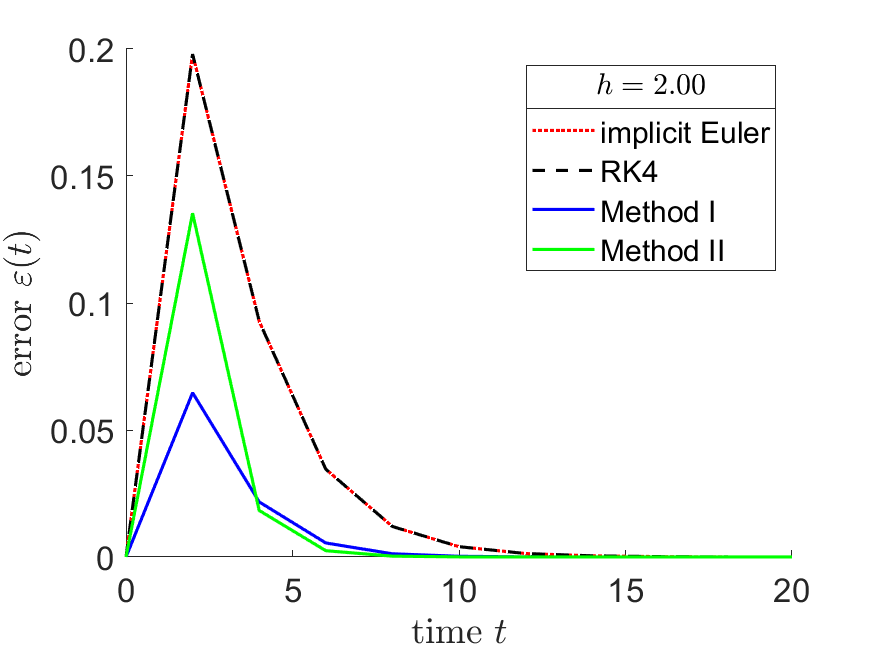}
\caption{}
\label{fig:fig1b}
\end{subfigure}
\begin{subfigure}{0.45\textwidth}
\centering
\includegraphics[width=\textwidth]{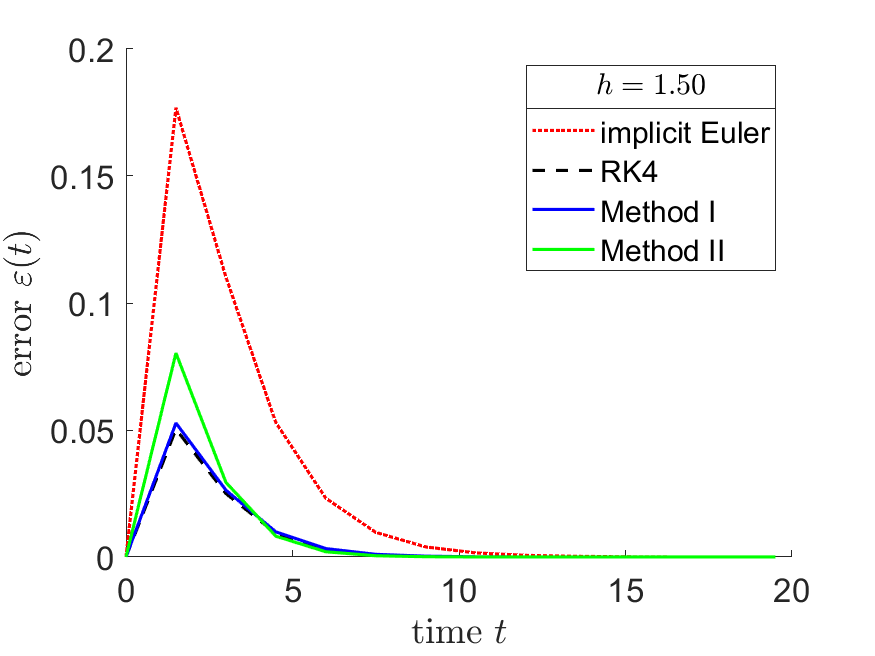}
\caption{}
\label{fig:fig1c}
\end{subfigure}
\begin{subfigure}{0.45\textwidth}
\centering
\includegraphics[width=\textwidth]{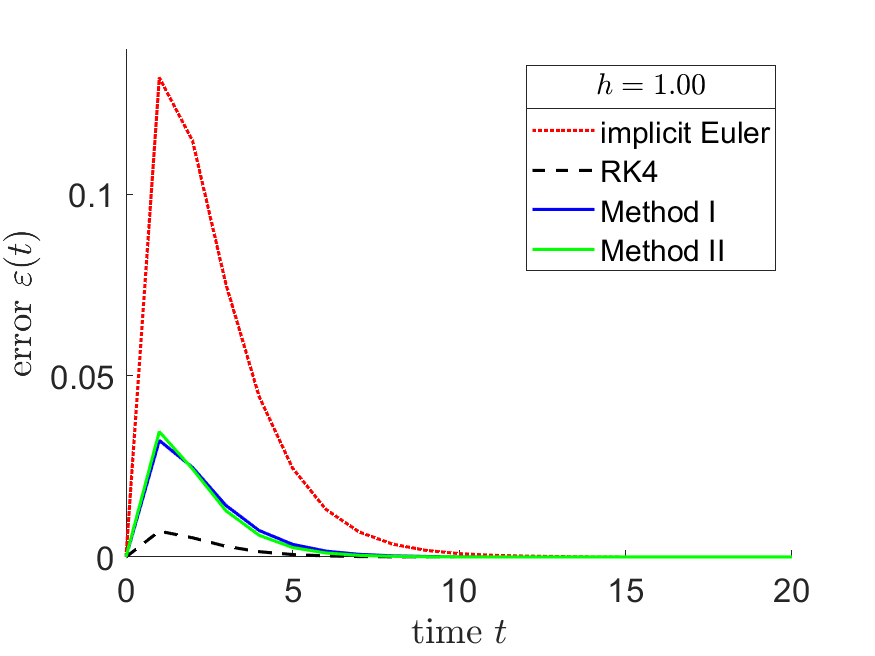}
\caption{}
\label{fig:fig1d}
\end{subfigure}
\begin{subfigure}{0.45\textwidth}
\centering
\includegraphics[width=\textwidth]{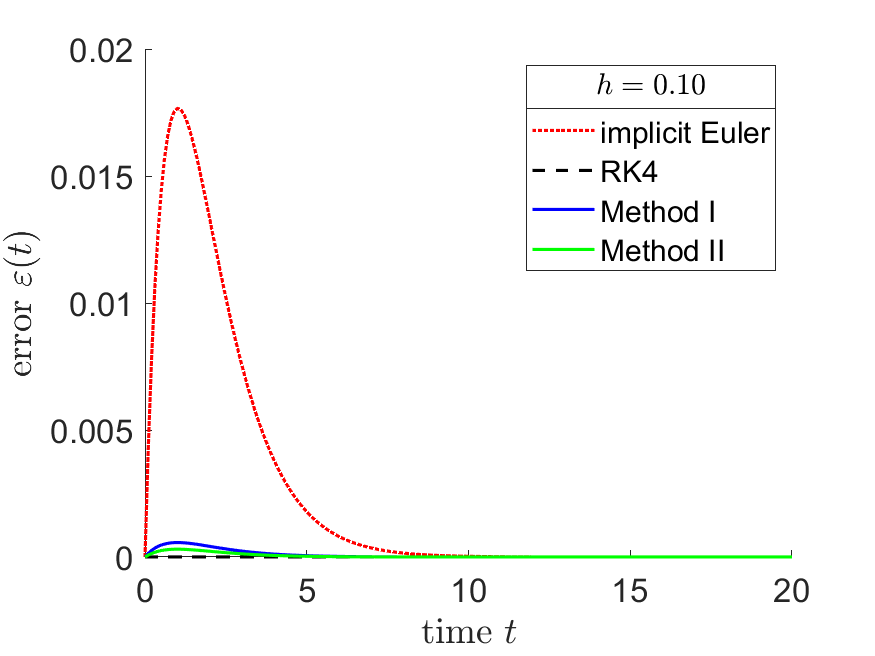}
\caption{}
\label{fig:fig1e}
\end{subfigure}
\begin{subfigure}{0.45\textwidth}
\centering
\includegraphics[width=\textwidth]{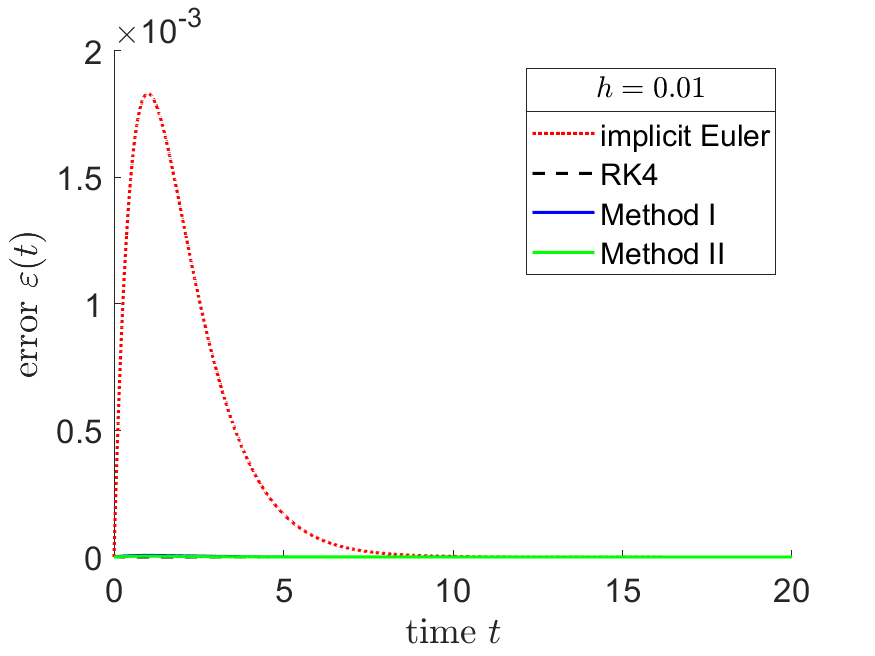}
\caption{}
\label{fig:fig1f}
\end{subfigure}
\caption{Representative numerical results for the test problem~\eqref{eq:test}, comparing the performance of Methods~I and~II against the implicit Euler method and the explicit fourth-order Runge-Kutta method (RK4) for \subref{fig:fig1a}~$h=2.50$, \subref{fig:fig1b}~$h=2.00$, \subref{fig:fig1c}~$h=1.50$, \subref{fig:fig1d}~$h=1.00$, \subref{fig:fig1e}~$h=0.10$, and \subref{fig:fig1f}~$h=0.01$. Here the error is the absolute difference between the numerical solution and the exact solution $e^{-t}$.}
\label{fig:test}
\end{figure}

\begin{figure}[!p]
\centering
\includegraphics[width=0.9\textwidth]{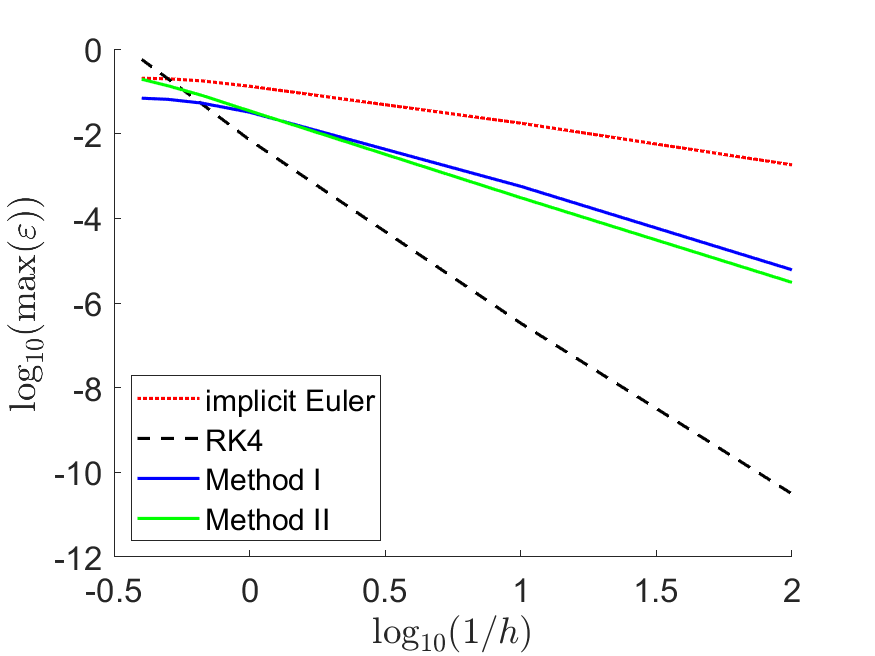}
\caption{Maximum absolute error versus inverse step size with logarithmic scaling for the results shown in Figure~\ref{fig:test}.}
\label{fig:test2}
\end{figure}

In the most general possible case, we could insist that any homogeneous function $\mathscr{F}$ of $p_{n}$ and $p_{n+1}$ vanishes, \emph{i.e.},
\begin{equation}
\mathscr{F}(p_{n},p_{n+1})=0.
\end{equation}
It is conceivable that a nonlinear constraint on $p_{n}$ and $p_{n+1}$ could result in improved performance, although again we leave such investigations for future work. For now we simply observe that, however we decide to constrain $p_{n}$ to $p_{n+1}$, by doing so we circumvent the need to specify an initial condition for $p_{0}$. That is, although the mathematically equivalent second-order formulation~\eqref{eq:secondorderEOM} generally requires an additional initial condition, a symplectic integrator for the associated Lagrangian~\eqref{eq:L} does not.

Figures~\ref{fig:test} and~\ref{fig:test2} show representative numerical results for the test problem~\eqref{eq:test}, comparing the performance of Methods~I and~II against the implicit Euler method and the explicit fourth-order Runge-Kutta method (RK4) for various time increments $h$. Here we are plotting the absolute error,
\begin{equation}
\varepsilon(t)\equiv |v(t) - e^{-t}|,
\end{equation}
between the numerical solution and the exact solution. In all cases, the error decays to zero as $t\rightarrow\infty$. It can be seen from Figures~\ref{fig:test} and~\ref{fig:test2} that, for sufficiently large time increments ($h\gtrsim 1.5$), Methods~I and~II outperform \emph{both} implicit Euler and RK4. For smaller time increments ($h\lesssim1.5$), Methods~I and~II still outperform implicit Euler by several orders of magnitude, but RK4 outperforms Methods~I and~II. Again, we expected Methods~I and~II to outperform implicit Euler, since these methods possess higher-order correction terms, and these results confirm our expectations. 

\section{Quadratic drag}\label{sec:quadraticdrag}

\begin{figure}[!p]
\centering
\begin{subfigure}{0.45\textwidth}
\centering
\includegraphics[width=\textwidth]{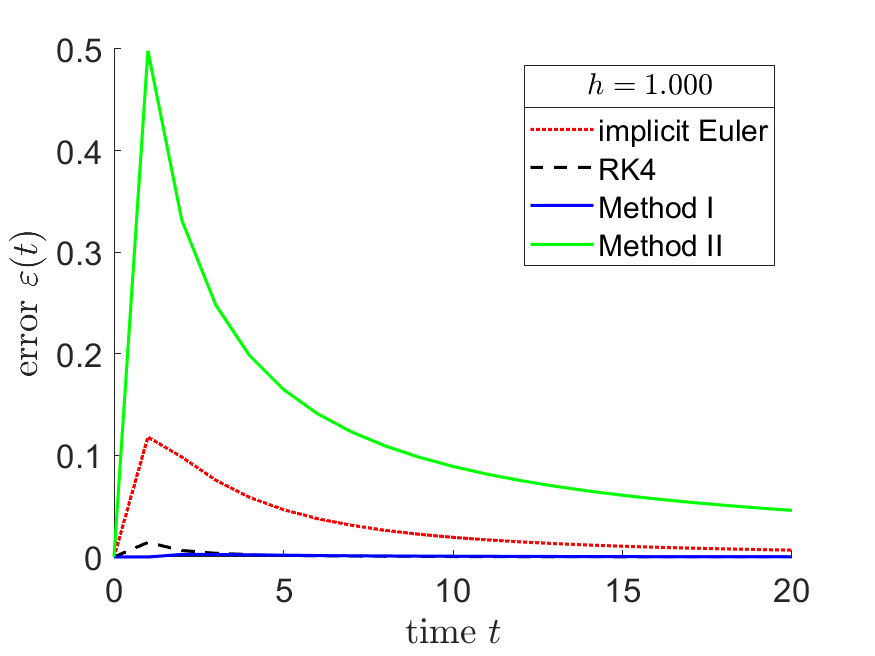}
\caption{}
\label{fig:fig3a}
\end{subfigure}
\begin{subfigure}{0.45\textwidth}
\centering
\includegraphics[width=\textwidth]{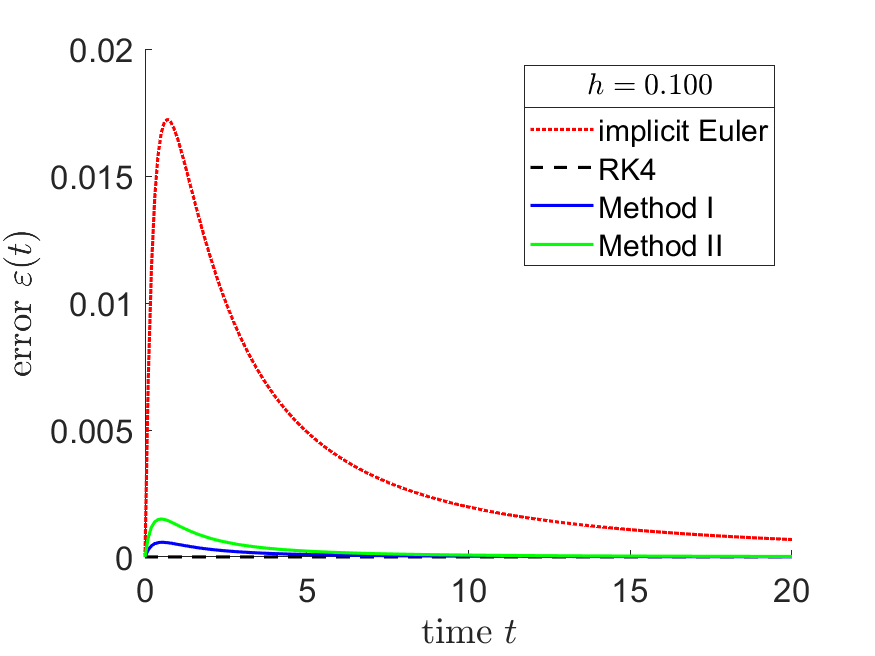}
\caption{}
\label{fig:fig3b}
\end{subfigure}
\begin{subfigure}{0.45\textwidth}
\centering
\includegraphics[width=\textwidth]{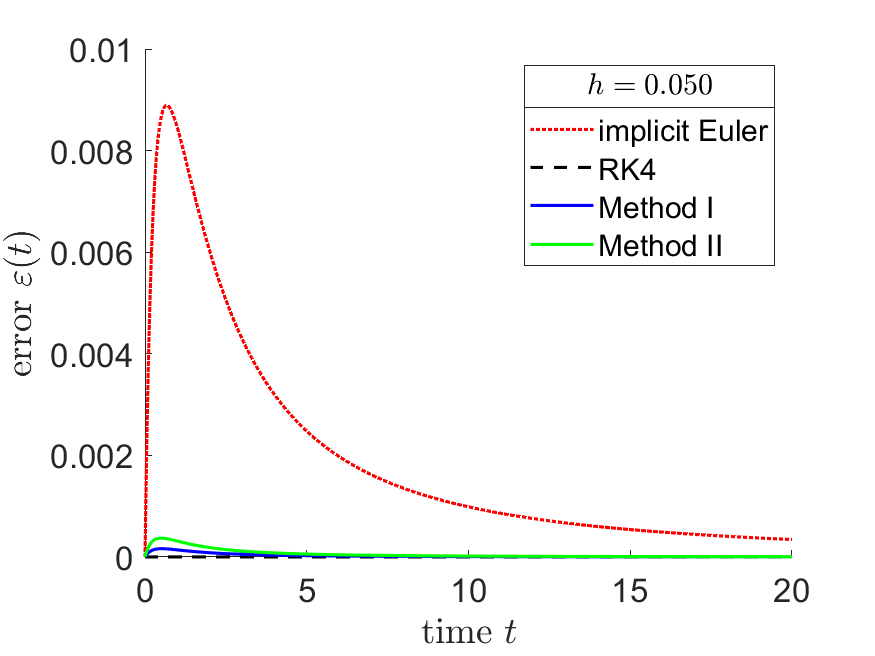}
\caption{}
\label{fig:fig3c}
\end{subfigure}
\begin{subfigure}{0.45\textwidth}
\centering
\includegraphics[width=\textwidth]{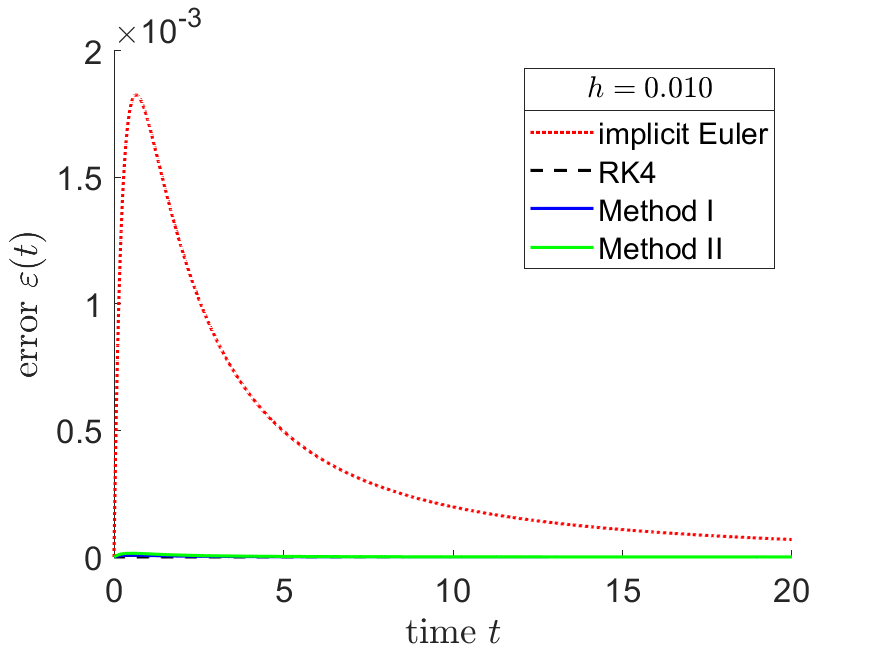}
\caption{}
\label{fig:fig3d}
\end{subfigure}
\begin{subfigure}{0.45\textwidth}
\centering
\includegraphics[width=\textwidth]{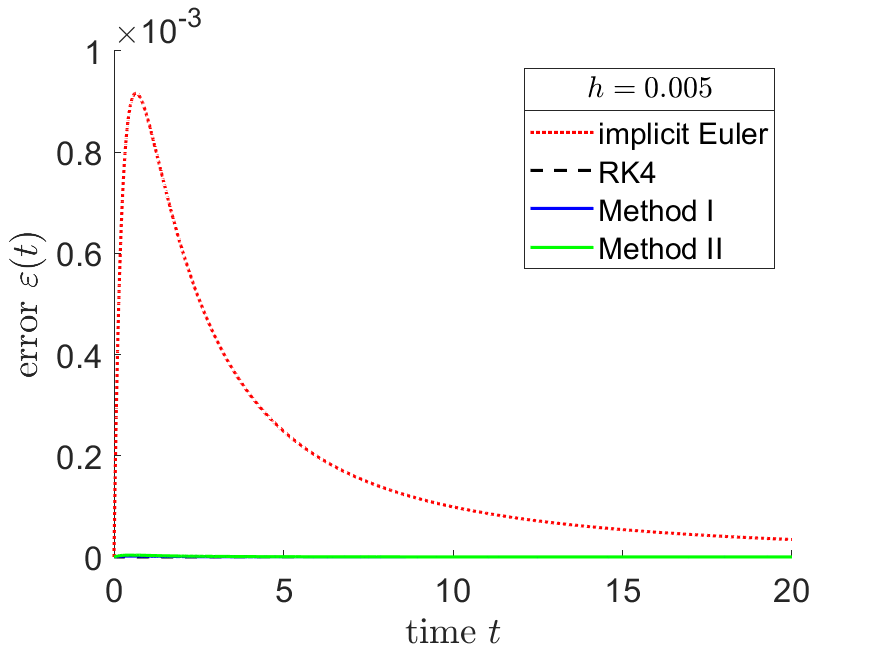}
\caption{}
\label{fig:fig3e}
\end{subfigure}
\begin{subfigure}{0.45\textwidth}
\centering
\includegraphics[width=\textwidth]{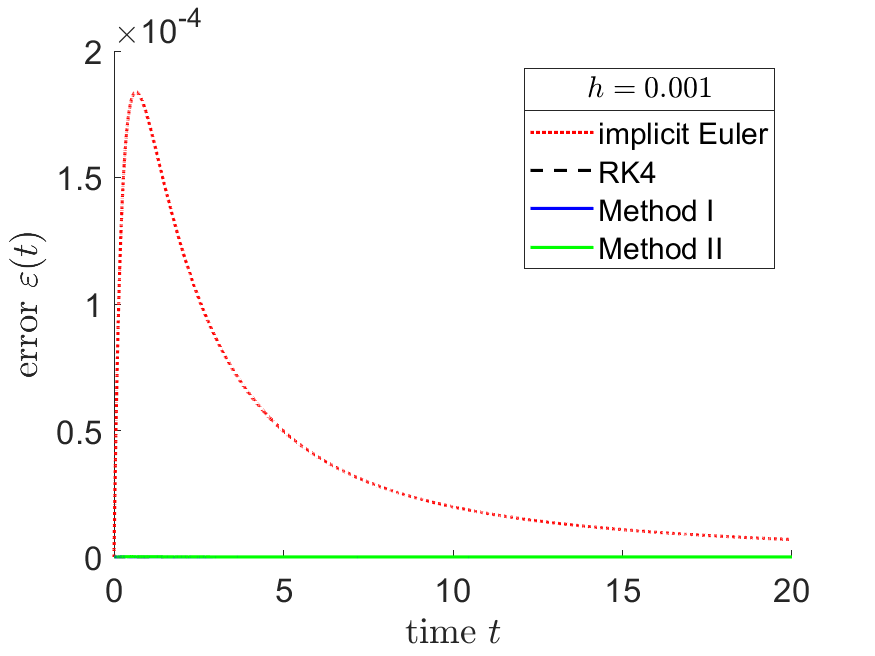}
\caption{}
\label{fig:fig3f}
\end{subfigure}
\caption{Representative numerical results for the nonlinear problem~\eqref{eq:quadraticdrag}, comparing the performance of Methods~I and~II against the implicit Euler method and the explicit fourth-order Runge-Kutta method (RK4) for \subref{fig:fig3a}~$h=1.000$, \subref{fig:fig3b}~$h=0.100$, \subref{fig:fig3c}~$h=0.050$, \subref{fig:fig3d}~$h=0.010$, \subref{fig:fig3e}~$h=0.005$, and \subref{fig:fig3f}~$h=0.001$. Here the error is the absolute difference between the numerical solution and the exact solution $(1+t)^{-1}$.}
\label{fig:quadraticdrag}
\end{figure}

\begin{figure}[!p]
\centering
\includegraphics[width=0.9\textwidth]{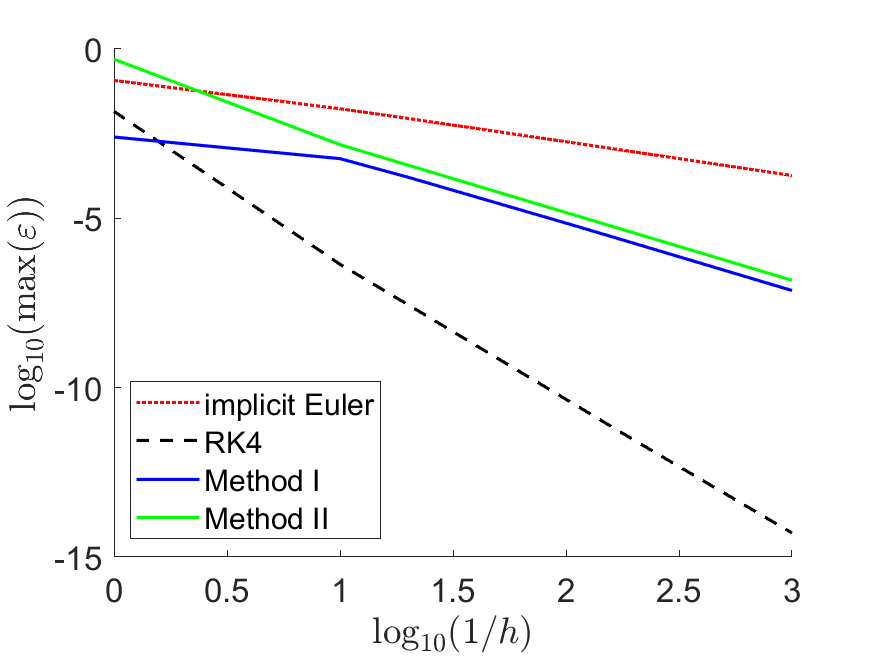}
\caption{Maximum absolute error versus inverse step size with logarithmic scaling for the results shown in Figure~\ref{fig:quadraticdrag}.}
\label{fig:quadraticdrag2}
\end{figure}

We now consider a nonlinear problem in the form of quadratic drag in one dimension:
\begin{equation}\label{eq:quadraticdrag}
\dot{v}=-v^{2}, \quad v(0)=1.
\end{equation}
This, too, is an intrinsically non-Hamiltonian problem to which symplectic integrators cannot be applied directly. However, proceeding as before, the Lagrangian of the mathematically equivalent second-order problem is
\begin{equation}
L(v,\dot{v})\equiv\frac{1}{2}\dot{v}^{2}+\dot{v}v^{2}+\frac{1}{2}v^{4}.
\end{equation}
To this Lagrangian we may apply the same Methods~I and~II from the previous section. 

Figures~\ref{fig:quadraticdrag} and~\ref{fig:quadraticdrag2} show analogous results to Figures~\ref{fig:test} and~\ref{fig:test2} for the nonlinear problem~\eqref{eq:quadraticdrag}. Here we are plotting the absolute error,
\begin{equation}
\varepsilon(t)\equiv \left|v(t) - \frac{1}{1+t}\right|,
\end{equation}
between the numerical solution and the exact solution. Once again, Methods~I and~II generally outperform implicit Euler, but for smaller time increments RK4 outperforms Methods~I and~II.

Before moving on, we pause again to emphasize that MacKay's method~\cite{MacKay1992} represents just one of the myriad existing symplectic integration schemes in the literature, and one of the simplest such schemes no less. We expect dissipative symplectic integrators based on more sophisticated schemes, such as that of Marsden~\&~West~\cite{Marsden2001}, to perform significantly better.

\section{Unsteady Poiseuille flow}\label{sec:poiseuille}

\begin{figure}[!b]
\centering
\includegraphics[width=0.5\textwidth]{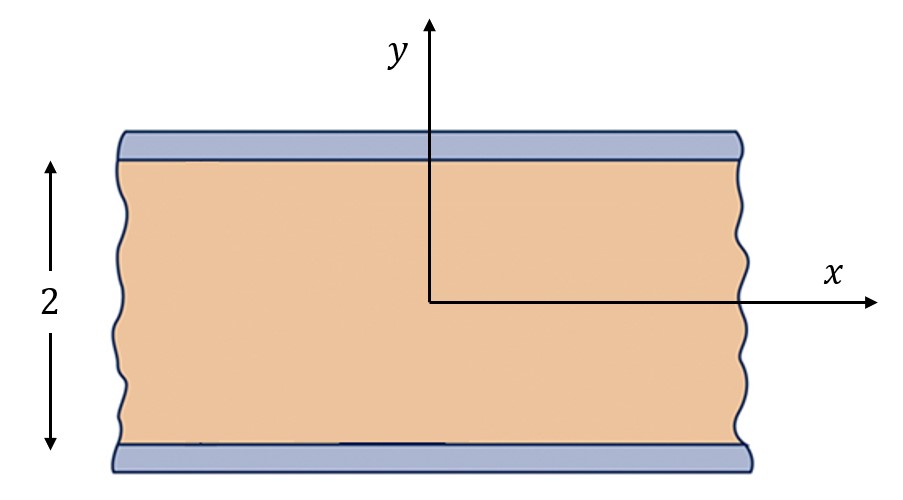}
\caption{Schematic diagram for closed channel flow, illustrating the coordinates and dimensions.}
\label{fig:channel}
\end{figure}

Like the simple dissipative problems considered in the previous sections, the Navier-Stokes equations~\cite{Stokes1845,LemarieRieusset2018} also lack symplectic structure. As such, symplectic integrators cannot be applied directly to viscous fluids modeled by the Navier-Stokes equations. Nevertheless, Sanders~\emph{et al.}~\cite{Sandersetal2024} recently applied the time-averaged principle of least squares to obtain a mathematically equivalent, canonical Hamiltonian formulation of the Navier-Stokes problem~\cite{Sanders2023b}. Notably, the higher-order formulation of Sanders~\emph{et al.}~\cite{Sandersetal2024} \emph{does} possess symplectic structure and so lends itself readily to symplectic integrators. To demonstrate the applicability of this approach to the numerical solution of the Navier-Stokes equations~\cite{Stokes1845,LemarieRieusset2018} and computational fluid dynamics in general, here we will apply the unconditionally stable symplectic integrators developed in the previous section (Methods~I and~II) to a simple viscous fluid flow problem. 

In particular, let us consider the case of developing viscous flow through a closed channel consisting of two infinite, flat plates, as shown schematically in Figure~\ref{fig:channel}. In what follows, we will non-dimensionalize all quantities using the fluid density $\rho$, fluid viscosity $\mu$, and channel half-width (so that the full channel width is 2, as shown). 

Let there be a fixed and known pressure gradient $\gamma$ along the length of the channel such that $p(x)=p_{0}-\gamma x$, where $p_{0}=p(0)$. Then the fluid will flow in the $x$-direction with velocity profile $u=u(y,t)$. Such a flow is automatically incompressible. The residual of the non-dimensionalized Navier-Stokes equation is
\begin{equation}\label{eq:NS}
\mathcal{R}\equiv-\gamma-\frac{\partial^{2}u}{\partial y^{2}}+\frac{\partial u}{\partial t}=0.
\end{equation}
Let the fluid be initially at rest, so that $u(y,0)=0$. If the plates are fixed, we have Poiseuille flow, and the boundary conditions are $u(\pm1,t)=0$. For reference, the exact solution is
\begin{equation}
u(y,t)=\sum_{n>0,\text{ odd}}\left(\frac{16\gamma}{n^{3}\pi^{3}}\right)\sin{\left(\frac{n\pi(y+1)}{2}\right)}\left[1-e^{-(n\pi/2)^{2}t}\right].
\end{equation}
The steady-state velocity profile, obtained by setting $\partial u/\partial t=0$ in~\eqref{eq:NS}, is easily found to be $u_{\infty}(y)=\tfrac{1}{2}\gamma(1-y^{2})$.

For the purposes of applying Methods~I and~II to this problem, here we will assume an approximate shape function satisfying the initial- and boundary conditions:
\begin{equation}\label{eq:shapefunction}
\hat{u}(y,t)=v(t)(1-y^{2}), \quad v(0)=0,
\end{equation}
where $v(t)=\hat{u}(0,t)$ is the maximum velocity at $y=0$. Again, this is only an approximate solution to the unsteady problem. Furthermore, although it has the same form as the steady-state solution, it would be a natural choice even without \emph{a priori} knowledge of the steady-state solution. {A more sophisticated approach would be to discretize the domain $y\in(-1,+1)$ into finite volumes and use shape functions to interpolate between the nodes. Again, for the present purposes, an approximate shape function over the entire domain will suffice to demonstrate the success of Methods~I and~II.} 

Substituting $\hat{u}(y,t)$ into the Navier-Stokes equation, we obtain
\begin{equation}\label{eq:Rpoiseuille}
\hat{\mathcal{R}}[v(t);y]\equiv-\gamma+2v(t)+\dot{v}(t)(1-y^{2})=0.
\end{equation}
The Lagrangian of the second-order formulation then takes the form
\begin{equation}
L(v,\dot{v})=\int_{-1}^{+1}\tfrac{1}{2}\hat{\mathcal{R}}^{2}\text{d}y=\gamma^{2}-4\gamma v+4v^{2}+\tfrac{4}{3}(-\gamma+2v)\dot{v}+\tfrac{8}{15}\dot{v}^{2}.
\end{equation}
From this we evaluate MacKay's~\cite{MacKay1992} discrete Lagrangian~\eqref{eq:discreteLMacKay}.

Employing Method~I (setting $p_{n+1}=0$), we obtain
\begin{equation}
v_{n+1}=\frac{\tfrac{8}{15}v_{n}+\gamma h^{2}+\tfrac{2}{3}\gamma h}{\tfrac{8}{15}+2h^{2}+\tfrac{4}{3}h}.
\end{equation}
We can see by inspection that this scheme is unconditionally stable, since in the limit as $h\rightarrow\infty$, we have
\begin{equation}
\lim_{h\rightarrow\infty}v_{n+1}=\lim_{h\rightarrow\infty}\frac{\tfrac{8}{15}v_{n}+\gamma h^{2}+\tfrac{2}{3}\gamma h}{\tfrac{8}{15}+2h^{2}+\tfrac{4}{3}h}=\lim_{h\rightarrow\infty}\frac{2\gamma h+\tfrac{2}{3}\gamma}{4h+\tfrac{4}{3}}=\lim_{h\rightarrow\infty}\frac{2\gamma}{4}=\frac{1}{2}\gamma,
\end{equation}
in agreement with the known steady-state velocity $u_{\infty}(0)=\tfrac{1}{2}\gamma$. Alternatively, employing Method~II (setting $p_{n+1}-p_{n}=0$), we obtain
\begin{equation}
v_{n+1}=\frac{(1-\tfrac{3}{2}h)v_{n}+\tfrac{3}{2}\gamma h}{1+\tfrac{3}{2}h}.
\end{equation}

\begin{figure}[!p]
\centering
\begin{subfigure}{0.45\textwidth}
\centering
\includegraphics[width=\textwidth]{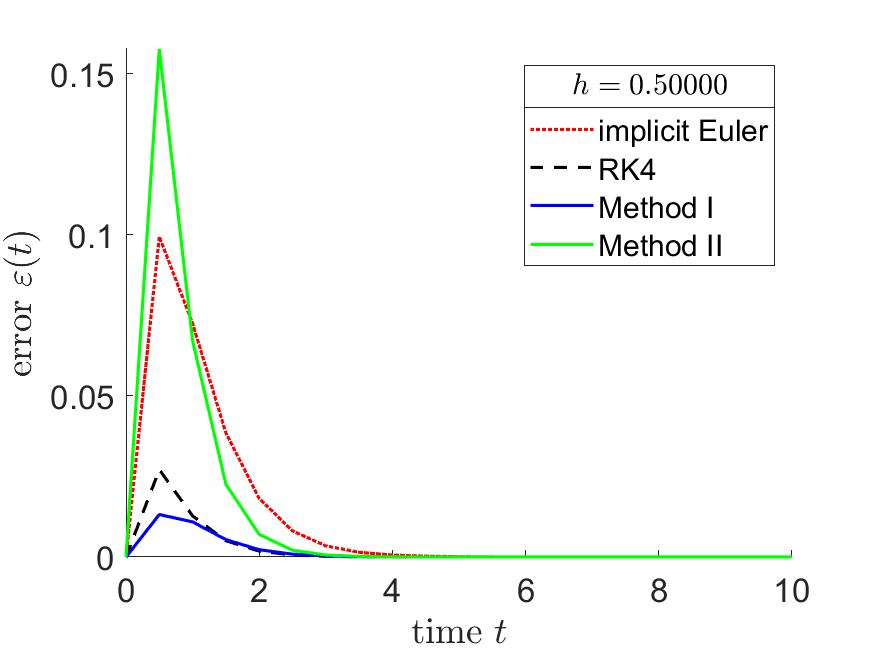}
\caption{}
\label{fig:fig6a}
\end{subfigure}
\begin{subfigure}{0.45\textwidth}
\centering
\includegraphics[width=\textwidth]{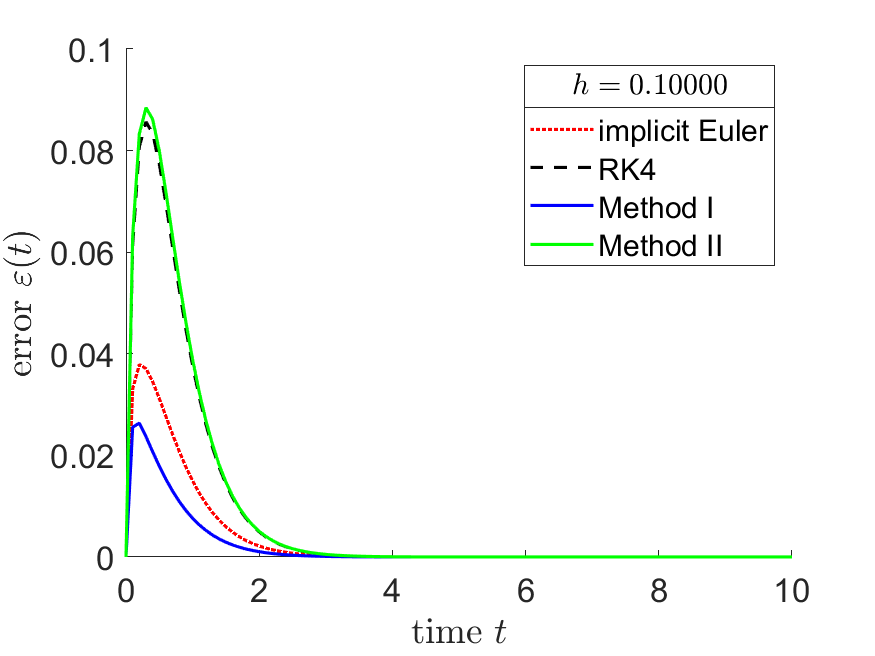}
\caption{}
\label{fig:fig6b}
\end{subfigure}
\begin{subfigure}{0.45\textwidth}
\centering
\includegraphics[width=\textwidth]{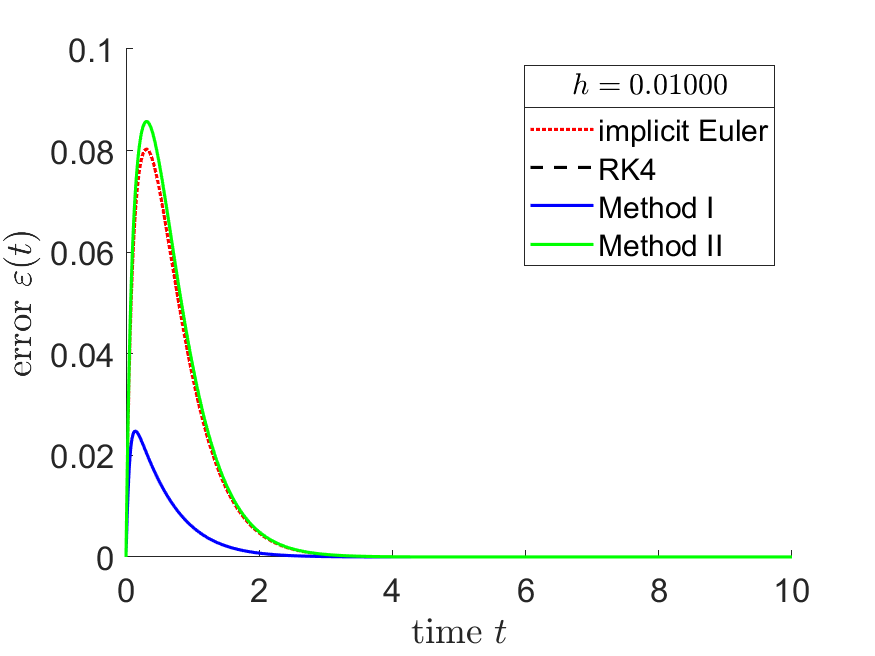}
\caption{}
\label{fig:fig6c}
\end{subfigure}
\begin{subfigure}{0.45\textwidth}
\centering
\includegraphics[width=\textwidth]{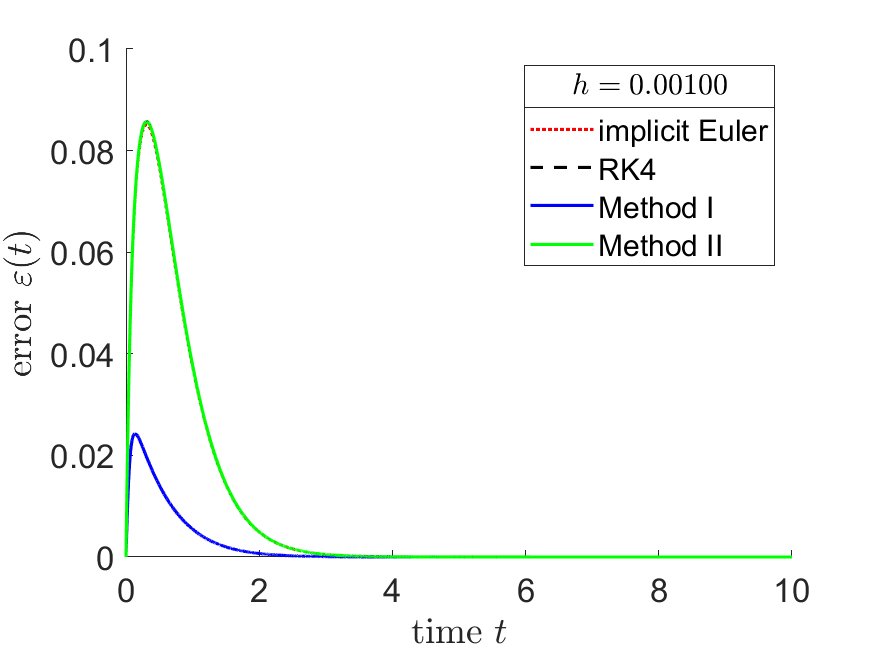}
\caption{}
\label{fig:fig6d}
\end{subfigure}
\begin{subfigure}{0.45\textwidth}
\centering
\includegraphics[width=\textwidth]{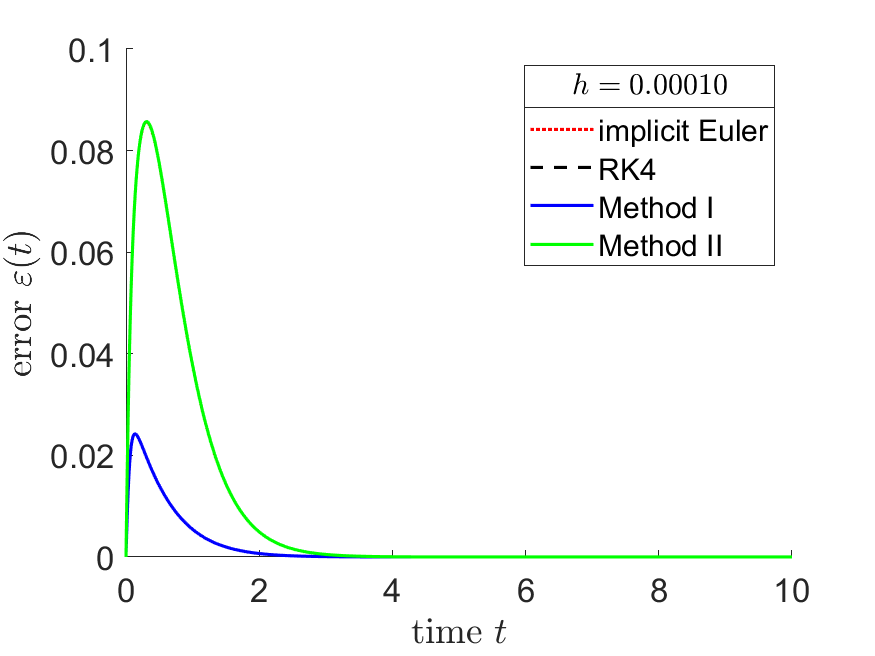}
\caption{}
\label{fig:fig6e}
\end{subfigure}
\begin{subfigure}{0.45\textwidth}
\centering
\includegraphics[width=\textwidth]{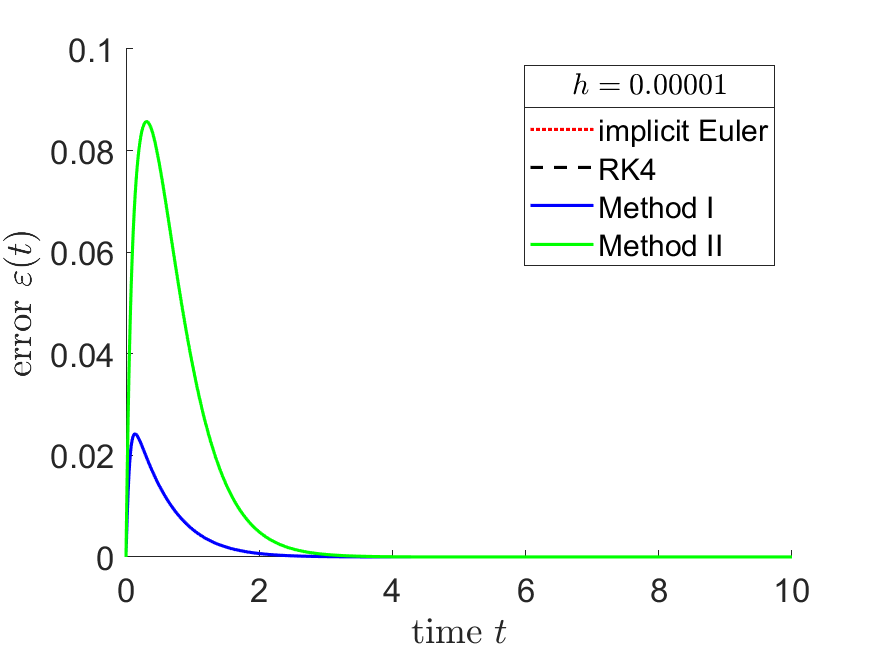}
\caption{}
\label{fig:fig6f}
\end{subfigure}
\caption{Representative numerical results for the unsteady Poiseuille flow problem~\eqref{eq:NS} with $\gamma=2$, comparing the performance of Methods~I and~II against the implicit Euler method and the explicit fourth-order Runge-Kutta method (RK4) for \subref{fig:fig6a}~$h=0.50000$, \subref{fig:fig6b}~$h=0.10000$, \subref{fig:fig6c}~$h=0.01000$, \subref{fig:fig6d}~$h=0.00100$, \subref{fig:fig6e}~$h=0.00010$, and \subref{fig:fig6f}~$h=0.00001$. Here the error is the absolute difference between the numerical solution and the exact solution for the maximum flow velocity at $y=0$.}
\label{fig:poiseuille}
\end{figure}

\begin{figure}[!p]
\centering
\includegraphics[width=0.9\textwidth]{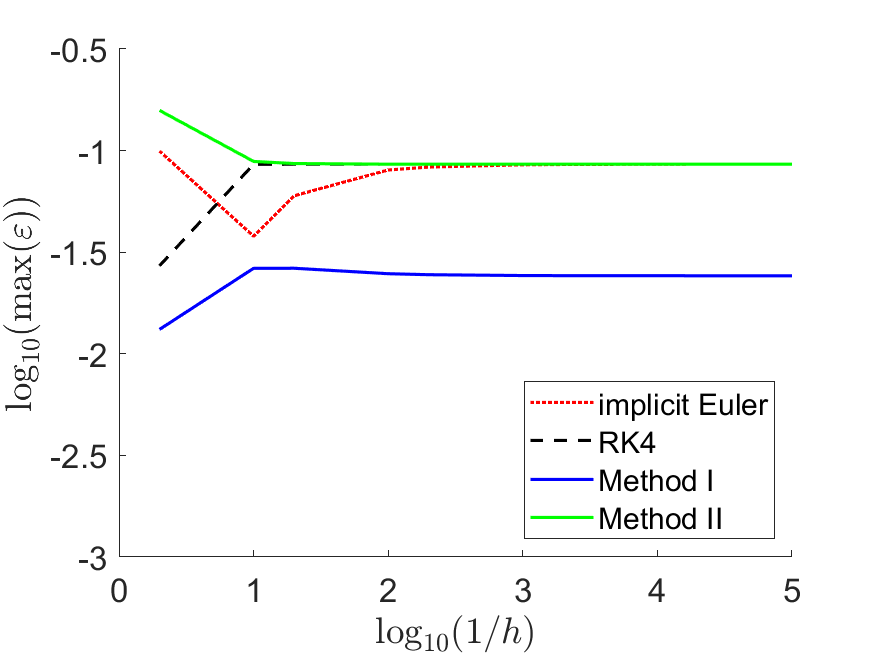}
\caption{Maximum absolute error versus inverse step size with logarithmic scaling for the results shown in Figure~\ref{fig:poiseuille}. The error does not converge to zero as the step size decreases because of the approximate form of the assumed shape function~\eqref{eq:shapefunction}. The numerical solution simply converges as closely as it can to the exact solution.}
\label{fig:poiseuille2}
\end{figure}

Figures~\ref{fig:poiseuille} and~\ref{fig:poiseuille2} show representative numerical results for $\gamma=2$ obtained using Methods~I and~II with various time increments $h$. For comparison, we consider solutions to the domain-averaged problem
\begin{equation}\label{eq:Rpoiseuille}
\int_{-1}^{+1}\hat{\mathcal{R}}\text{d}y=-2\gamma+4v(t)+\tfrac{4}{3}\dot{v}(t)=0,
\end{equation}
obtained using the implicit Euler and RK4 methods. Here we are plotting the absolute error between the numerical solution $v(t)$ and the exact velocity $u(0,t)$ at $y=0$:
\begin{equation}
\varepsilon(t)\equiv \left|v(t)-\sum_{n>0,\text{ odd}}\left(\frac{16\gamma}{n^{3}\pi^{3}}\right)\sin{\left(\frac{n\pi}{2}\right)}\left[1-e^{-(n\pi/2)^{2}t}\right]\right|.
\end{equation}
In all cases, the error decays to zero as $t\rightarrow\infty$. Due to the approximate form of the shape function, the error does not converge to zero for small times. Instead, the numerical solution converges as closely as it can to the exact solution. Note that, with $\gamma=2$, $u_{\infty}(0)=1$, so that the absolute error coincides with the relative error for sufficiently large $t$.

Interestingly---and perhaps surprisingly---this time Method~I generally outperforms \emph{all three other methods}, even for very small values of $h$. With $h\lesssim0.010$, the maximum error exhibited by Method~I is only roughly one-quarter that of the other three methods. To the authors' knowledge, this is the very first time that a symplectic integrator has been applied successfully to the Navier-Stokes equations.

The success of the present framework at constructing accurate and unconditionally stable symplectic integrators for the Navier-Stokes equations provides direct empirical validation of the canonical Hamiltonian formulation of the Navier-Stokes problem recently published by Sanders~\emph{et al.}~\cite{Sandersetal2024}. The present authors also see in these results great potential for improved numerical schemes for computational fluid dynamics. Once again, we emphasize that the present results were obtained using one of the simplest symplectic integrators in the literature~\cite{MacKay1992}. We expect more sophisticated symplectic integration schemes~\cite{Marsden2001} to perform significantly better, and we imagine that such schemes may one day be incorporated into commercial computational fluid dynamics programs. 

\section{Conclusion}\label{sec:conclusion}

This paper has presented a general framework for applying existing symplectic integration schemes \cite{Hairer2006,Maeda1980,Maeda1982,Veselov1988,Veselov1991,Suris1988,Suris1990,Suris1996,Suris2003,MacKay1992,Wendlandt1997,Kane2000,Marsden2001,Lew2003,Lew2004,Lew2004a,Kale2007,OberBloebaum2011,Lew2016} to dynamical systems that lack symplectic structure, with particular emphasis on the Navier-Stokes equations~\cite{Stokes1845,LemarieRieusset2018,Sandersetal2024}. Two such schemes, obtained by applying the method of MacKay~\cite{MacKay1992} to the Lagrangian~\eqref{eq:L}, have been shown to be unconditionally stable. Additionally, for the unsteady flow problem considered here, one of these methods (Method~I) outperforms both the implicit Euler method and the fourth-order Runge-Kutta method in terms of accuracy for a given time step. If MacKay's method~\cite{MacKay1992} can be made to outperform both implicit Euler and RK4 for dissipative systems, we expect unprecedented performance from more sophisticated symplectic integrators, such as that of Marsden~\&~West~\cite{Marsden2001}. Future work will focus on the development of such schemes for the Navier-Stokes equations~\cite{Stokes1845,LemarieRieusset2018,Sandersetal2024}. Perhaps one day such schemes may be incorporated into commercial computational fluid dynamics programs. Furthermore, we interpret the numerical success of the present method as direct empirical validation of the canonical Hamiltonian formulation of the Navier-Stokes problem published by Sanders~\emph{et al.}~\cite{Sandersetal2024}.

%
%
%
%
%
%

%
%
%
%

\begin{spacing}{0.25}
\small
\bibliographystyle{asmems4}
\bibliography{Sandersetal2024arXiv}
\end{spacing}

\end{document}